\documentclass[12pt]{amsart}

\usepackage[dvips]{graphicx}
	\usepackage{amsmath,amssymb,amsthm,wrapfig,amsfonts,enumerate,latexsym}

	\newtheorem{dfn}{Definition}[section]
	\newtheorem{thm}[dfn]{Theorem}
	\newtheorem{prop}[dfn]{Proposition}
	\newtheorem{cor}[dfn]{Corollary}
	\newtheorem{lem}[dfn]{Lemma}
	
	\newtheorem{ex}[dfn]{Example}

	\newtheorem{ack}{Acknowledgements\!\!}

	\newcounter{yon}
	\numberwithin{equation}{section}

	\def\notin{\not\in}

	\newcommand{\dist}{\mathop{\mathit{d}} \nolimits}

	\newcommand{\grad}{\mathop{\mathrm{grad}}                   \nolimits}

	\newcommand{\obin}{\mathop{\mathrm{Obs}L^p\mathrm{\text{-}Var}}
	\nolimits}
	\newcommand{\obinin}{\mathop{\mathrm{Obs}L^2\mathrm{\text{-}Var}}
	\nolimits}
	
	\newcommand{\obvar}{\mathop{\mathrm{Obs}L^{1}\mathrm{\text{-}Var}}
	\nolimits}

\begin{document}

	\title[An asymptotic variant of the Fubini theorem for maps into CAT(0)-spaces]
    {An asymptotic variant of the Fubini theorem\\ for maps into CAT(0)-spaces}
	\author[Kei Funano]{Kei Funano}
	\address{Mathematical Institute, Tohoku University, Sendai 980-8578, JAPAN}
	\email{sa4m23@math.tohoku.ac.jp}
	\thanks{This work was partially supported by Research Fellowships of
	the Japan Society for the Promotion of Science for Young Scientists.}
	\dedicatory{}
	\date{\today}

	\maketitle


\begin{abstract}The classical Fubini theorem asserts that the multiple
 integral is equal to the repeated one for any integrable function on a
 product measure space. In this paper, we derive
    an asymptotic variant of the Fubini theorem for maps into CAT(0)-spaces from
    the $L^1$ and $L^2$-concentration of the maps.
 \end{abstract}
	\setlength{\baselineskip}{5mm}

    \section{Introduction and statement of the main result}

The classical Fubini theorem asserts that the multiple integral is equal
to the repeated one for any integrable function on a product measure
space. In this paper, we prove an asymptotic variant of the Fubini theorem for maps into CAT(0)-spaces.

For this purpose, let us define the expectation (integral) for a map from
 a probability space into a CAT(0)-space. Throughout this section, let $N$ be a CAT(0)-space. Let $(\Omega ,
 \mathcal{A}, \mathbb{P})$ be a Probability space. For an
 $N$-valued random variable $Z:\Omega \to N$ such that the push-forward
 measure $Z_{\ast}\mathbb{P}$ has finite moment of order $1$, we
 define its \emph{expectation} $\mathbb{E}(Z)$ as the center of mass of the measure
 $Z_{\ast}\mathbb{P}$ (see Subsection 2.2 for the definition of the
 center of mass). This definition of expectations is based on the
 classical point of view of \cite{gaus}. In \cite{gaus}, C. F. Gauss
 defined the expectations of random variables with values in Euclidean spaces as
 the above way. In the context of metric spaces,
 this point of view was successfully used by \cite{fre}, \cite{jost2},
 \cite{sturm}, and many others. 

Let $(X,\dist_X,\mu_X)$ and $(Y,\dist_Y,\mu_Y)$ be two mm-spaces. Here, an
    \emph{mm-space} is a triple
    $(X,\dist_X,\mu_X)$ of a set $X$, a complete separable
    distance function $\dist_X$ on $X$, and a Borel probability measure
    $\mu_X$ on $(X,\dist_X)$ with full-support. The \emph{product
    mm-space} of $X$ and $Y$ is the mm-space $X\times Y$ equipped with the $\ell_2$-distance
    function and the product probability measure. For a
 Borel measurable map $f:X\times Y \to N$ such that the push-forward
measure $f_{\ast}(\mu_X \times \mu_Y)$ has finite
moment of order $1$ and $y \in Y$, we shall consider the map
$f^y:X \to N$ defined by
$f^{y}(x):=f(x,y)$. Note that the
push-forward measure $(f^{y})_{\ast}(\mu_X)$ has finite moment of order $1$ for $\mu_Y$-a.e. $y \in Y$.
For defining the repeated integral for the map $f$, we assume the following:

(1) The map $g_f :Y \to N$ defined by $g_f(y):= \mathbb{E}(f^y)$ is
    Borel measurable.

(2) The push-forward measure $(g_f)_{\ast}(\mu_Y)$
    have finite moment of order $1$.

    If the map $f$ is uniformly continuous, then
    the map $g_f$ satisfies the above $(1)$ and $(2)$ (see Lemmas
    \ref{l3} and \ref{l4}).
    It seems that the above $(1)$ and $(2)$ hold for an arbitrary Borel
    measurable map $f$, but the author does not know how to prove it as
    of now. For a Borel measurable map $f$ satisfying the above (1) and (2), we define its \emph{repeated integral}
 $\mathbb{E}_y(\mathbb{E}(f^y))$ by the expectation
 $\mathbb{E}(g_f)$. We will see that the Fubini theorem
 $\mathbb{E}(f) = \mathbb{E}_y (\mathbb{E}(f^y))$ does not hold in
 general for a nonlinear CAT(0)-space $N$ (see Example \ref{e1}). However, we
 succeed to estimate their difference $\dist_N(\mathbb{E}(f) ,\mathbb{E}_y
 (\mathbb{E}(f^y)))$ by the term of the $L^1$ and $L^2$-variation of the map $f$:  Let $X$ be
 an mm-space and $p\geq 1$. Given a Borel measurable map $f:X\to N$, we define its
 \emph{$L^p$-variation} by 
\begin{align*}
V_p(f):=  \Big(\int \int_{X\times X}
 \dist_N(f(x),f(x'))^p    d\mu_X(x)d \mu_X(x')\Big)^{1/p}.
\end{align*}A main theorem of this paper is the following:

    \begin{thm}\label{th1}
     Let $X$ and $Y$ be two mm-spaces. Then, for any uniformly
     continuous map $f:X\times Y \to N$ such that the push-forward
     measure $f_{\ast}(\mu_X \times \mu_Y)$ has the finite moment of
     order $1$, we have
     \begin{align}\label{s1}
      \dist_N(\mathbb{E}(f), \mathbb{E}_y(\mathbb{E}(f^y))) \leq V_1(f)
      \end{align}
     and
     \begin{align}\label{s2}
      \dist_N(\mathbb{E}(f), \mathbb{E}_y(\mathbb{E}(f^y)))
   \leq \frac{1}{\sqrt{3}}V_2(f)
      \end{align}
     \end{thm}Jensen's inequality easily leads to the inequality
   (\ref{s1}). In the proof of the inequality (\ref{s2}), we
     iterate some K-T. Sturm's inequality about the center of mass of a
     probability measure on a CAT(0)-space (see Proposition \ref{p3}). We
     emphasize that the coefficient $1/\sqrt{3}$ of the inequality (\ref{s2}) cannot be obtained only
     from the inequalities (\ref{s1}) and $V_1(f)\leq V_2(f)$.

 Let $\{ X_n\}_{n=1}^{\infty}$ be a sequence of mm-spaces and
 $\{N_n\}_{n=1}^{\infty}$ a sequence of CAT(0)-spaces. For $p\geq 1$, we say that a sequence $\{ f_n:X_n \to N_n\}_{n=1}^{\infty}$ of Borel
 measurable maps \emph{$L^p$-concentrates} if $V_p(f_n )\to 0$ as $n\to
 \infty$. From the inequality (\ref{s1}), the $L^p$-concentration of uniformly
 continuous maps implies that the Fubini theorem for the maps
 ``almostly'' holds. The $L^2$-concentration theory of maps into CAT(0)-spaces was first
    studied by M. Gromov in \cite{gromovcat}. In \cite{funano2}, the author also studied
    relationships between the L\'{e}vy-Milman concentration
    theory of $1$-Lipschitz maps and the $L^p$-concentration theory of
    $1$-Lipschitz maps (see \cite{ledoux}, \cite{mil}, \cite{mil2},
 \cite{milsch}, and \cite{sch} for further information about the
 L\'{e}vy-Milman concentration theory).
Motivated by Gromov's works in \cite{gromovcat}, \cite{gromov2}, and
 \cite{gromov}, the
    author studied the $L^p$-concentration theory of $1$-Lipschitz maps into
    Hadamard manifolds and $\mathbb{R}$-trees in \cite{funano2} and
    \cite{funano1}. Combining Theorem \ref{th1} with author's works and
 Gromov's works, we obtain the following corollary: We shall consider each compact connected Riemannian manifold $M$ as an
    mm-space equipped with the volume measure normalized to have the total
    volume $1$. We denote by $\lambda_1(M)$ the non-zero first
 eigenvalue of the Laplacian on $M$.
    \begin{cor}\label{c1}Let $M$ be a compact Riemannian manifold. Then, for any
     $n$-dimensional Hadamard manifold $N'$ and $1$-Lipschitz map $f:M\times M\to N'$, we have
     \begin{align*}
      \dist_{N'}(\mathbb{E}(f), \mathbb{E}_y(\mathbb{E}(f^y))) \leq
      2 \sqrt{\frac{2n}{3\lambda_1(M)}}.
      \end{align*}For an $\mathbb{R}$-tree $T$ and a $1$-Lipschitz map
     $f:M\times M \to T$, we also have
     \begin{align*}
      \dist_T(\mathbb{E}(f), \mathbb{E}_y(\mathbb{E}(f^y)))^2 \leq \frac{8(38+16\sqrt{2})}{3\lambda_1(M)}. 
      \end{align*}
     \end{cor}

     \section{Preliminaries}

 \subsection{The Wasserstein distance function of order 1}
    Let $(X,\dist_X)$ be a complete metric space. For
$p\geq 1$, we
indicate by $\mathcal{P}_p(X)$ the set of all probability measures $\nu$ such
that $\nu$ has the separable support and $\int_{X}\dist_X(x,y)^p \ d\nu (y)<+\infty$ for
some (hence all) $x\in X$.

For $\mu,\nu \in \mathcal{P}_1(X)$, we
  define the \emph{Wasserstein distance $\dist_1^W(\mu,\nu)$ of order $1$} between $\mu$
  and $\nu$ as the infimum of $\int_{X\times X}\dist_X(x,y) \ d\pi
  (x,y)$, where $\pi \in \mathcal{P}_1(X\times X)$ runs over all \emph{couplings}
  of $\mu$ and $\nu$, that is, the probability measures $\pi$ with the property that $\pi (A\times
  X)=\mu(A)$ and $\pi (X\times A)=\nu (A)$ for any Borel subset
  $A\subseteq X$.

  \begin{thm}[{L. V. Kantorovich, cf.~\cite[Theorem
   5.1, Remark 6.5]{villani}}]\label{t1}For any $\mu,\nu\in \mathcal{P}_1(X)$, we have
   \begin{align*}
    \dist_{1}^W(\mu,\nu)= \sup \Big\{ \int_X \psi(x)  d\mu
    (x)-\int_X\psi (x) d\nu (x)\Big\}, 
    \end{align*}where the supremum is taken over all $1$-Lipschitz
   function $\psi:X\to \mathbb{R}$.
   \end{thm}

  \subsection{Basics of the center of mass of a measure on CAT(0)-spaces}

In this subsection, we review Sturm's works about probability measures on a
CAT(0)-spaces, which is needed for the proof of the main theorem. Refer \cite{jost} and \cite{sturm} for details.

We shall recall
    some standard terminologies in metric geometry. Let $(X,\dist_X)$ be
    a metric space. A rectifiable curve $\gamma:[0,1]\to X$ is called a
    \emph{geodesic} if its arclength coincides with
    the distance $\dist_X(\gamma(0),\gamma(1))$ and it has a constant speed,
    i.e., parameterized proportionally to the arc length. We say that
    $(X,\dist_X)$ is a \emph{geodesic metric space} if any two points in $X$ are joined
    by a geodesic between them. A geodesic metric
space $N$ is called a \emph{CAT(0)-space} if we have
\begin{align*}
\dist_N(x,\gamma (1/2))^2 \leq \frac{1}{2}\dist_N(x,y)^2 +
 \frac{1}{2}\dist_N(x,z)^2 - \frac{1}{4} \dist_N(y,z)^2
\end{align*}for any $x,y,z \in N$ and any minimizing geodesic $\gamma
 :[0,1]\to N$ from $y$ to $z$. For example, Hadamard manifolds, Hilbert
 spaces, and $\mathbb{R}$-trees are all CAT(0)-spaces. 

For any $\nu \in \mathcal{P}_1(X)$ and $z\in X$, we consider the function
$h_{z,\nu}:X\to \mathbb{R}$ defined by
\begin{align*}
h_{z,\nu}(x):= \int_{X} \{\dist_X(x,y)^2-\dist_X(z,y)^2 \}  d\nu(y).
\end{align*}
Note that
\begin{align*}
\int_X |\dist_X(x,y)^2-\dist_X(z,y)^2|  d\nu(y)\leq \dist_X(x,z)\int_X
 \{\dist_X(x,y)+ \dist_X(z,y)\} d \nu(y) <+\infty.
\end{align*}A point $z_0 \in X$ is called the
 \emph{center of mass} of the measure $\nu \in \mathcal{B}_1(X)$ if for
 any $z\in X$, $z_0$ is a unique minimizing
 point of the function $h_{z,\nu}$. We denote the point $z_0$ by
 $c(\nu)$. Note that if the measure $\nu $ moreover satisfies that
 $\nu \in \mathcal{P}_2(X)$, then we have
 \begin{align*}
  \int_X \dist_X(c(\nu),y)^2 \ d\nu(y) = \inf_{x\in X}\int_X
  \dist_X(x,y)^2 \ d\nu(y).
  \end{align*}
A metric space $X$ is said to be \emph{centric} if every $\nu \in
 \mathcal{P}_1(X)$ has the center of mass. 

\begin{prop}[{cf.~\cite[Proposition
 $4.3$]{sturm}}]\label{p1}
A CAT(0)-space is centric.
\end{prop}

A simple variational argument implies the following lemma:

\begin{lem}[{cf.~\cite[Proposition 5.4]{sturm}}]\label{l1}Let $H$ be a Hilbert space. Then, for each $\nu \in
 \mathcal{P}_1(H)$, we have
 \begin{align*}
  c (\nu)=  \int_H y  d\nu(y).
  \end{align*}
\end{lem}

\begin{lem}[{cf.~\cite[Proposition 5.10]{sturm}}]\label{l2}Let $N$ be a Hadamard
 manifold and $\nu \in \mathcal{P}_1(N)$. Then, $x= c(\nu)$ if and only
 if
 \begin{align*}
  \int_N \exp_x^{-1}(y)d\nu(y)=0.
  \end{align*}In particular, identifying the tangent space of $N$ at
 $c(\nu)$ with the Euclidean space of the same dimension, we have $c((\exp_{c(\nu)}^{-1})_{\ast}(\nu))=0$.
 \end{lem}
Let $(\Omega, \mathcal{A}, \mathbb{P})$ be a probability space and $N$
a centric metric space. For an $N$-valued random variable
$Z:\Omega \to N$ satisfying $Z_{\ast}\mathbb{P}\in \mathcal{P}_1(N)$, we
define its \emph{expectation} $\mathbb{E}(Z)\in N$ by the point
$c(Z_{\ast}\mathbb{P})$.

 Let $X$ be a geodesic metric space. A function $\varphi:X\to \mathbb{R}$
  is called \emph{convex} if the function $\varphi \circ \gamma
  :[0,1]\to \mathbb{R}$ is convex for each geodesic $\gamma : [0,1]\to X$.
\begin{prop}[{Convexity of a distance function, cf.~\cite[Corollary
   2.5]{sturm}}]\label{p2}Let $N$ be a CAT(0)-space and $\gamma, \eta:[0,1]\to N$
   be two geodesics. Then, for any $t\in [0,1]$, we have
   \begin{align*}
    \dist_N(\gamma(t),\eta(t)) \leq (1-t)\dist_N(\gamma(0), \eta(0)) + t
    \dist_N(\gamma (1),\eta(1)).
    \end{align*}
   \end{prop}

  \begin{thm}[{Jensen's inequality, cf.~\cite[Theorem
   6.2]{sturm}}]\label{t2}Let $N$ be a CAT(0)-space. Then, for any lower semicontinuous convex function $\varphi : N\to \mathbb{R}$ and $\nu\in \mathcal{P}_1(N)$, we have
 \begin{align*}
  \varphi (c(\nu))\leq \int_N \varphi (x) \ d\nu (x),
  \end{align*}provided the right-hand side is well-defined.
   \end{thm}

   Applying Proposition \ref{p2} to Theorem \ref{t2}, we obtain the following corollary:
   \begin{cor}\label{c2}Let $N$ be a CAT(0)-space. Then, for any $p_0 \in N$ and
    $\nu \in \mathcal{P}_1(N)$, we have
    \begin{align*}
     \dist_N(p_0,c(\nu)) \leq \int_N \dist_N(p_0, p)d\nu(p).
     \end{align*}
    \end{cor}

\begin{prop}[{Variance inequality, \cite[Proposition 4.4]{sturm}}]\label{p3}Let
 $N$ be a CAT(0)-space and $\nu \in \mathcal{P}_1(N)$. Then, for any
 $z\in N$, we have
 \begin{align}\label{s3}
  \int_N\{\dist_N(z,x)^2-\dist_N(c(\nu),x)^2\}  d\nu(x)
  \geq \dist_N(z,c(\nu))^2.
  \end{align}
 \end{prop}
Note that if $N$ is a Hilbert space, then we have the equality in (\ref{s3}).
 \begin{prop}[{cf.~\cite[Theorem 2.5]{sturm}}]\label{p4}Let $N$ be a
  CAT(0)-space. Then, for any $\mu, \nu \in
  \mathcal{P}_1(N)$, we have $\dist_N(c(\mu),c(\nu)) \leq \dist^W_1(\mu,\nu)$.
  \end{prop}

 \section{Proof of the main theorem}

 Let $X$ and $Y$ be a two mm-spaces and $N$ a CAT(0)-space.
 Given a uniformly continuous map $f:X\times Y \to N$ with $f_{\ast}(\mu_X \times \mu_Y)\in \mathcal{P}_1(N)$, we easily see
 that $(f^{y})_{\ast}(\mu_X) \in \mathcal{P}_1(N)$ for
 $\mu_Y$-a.e. $y\in Y$. Since $Y$ has the full-support and the map $f$
 is uniformly continuous, we see that $(f^{y})_{\ast}(\mu_X) \in
 \mathcal{P}_1(N)$ for any $y\in N$. We shall consider the map $g_f:Y\to N$ defined
 by $g_f(y):= \mathbb{E}(f^y)$.
 \begin{lem}\label{l3}The map $g_f: Y \to  N$ is uniformly
  continuous. In particular, the map is Borel measurable.
  \begin{proof}From Theorem \ref{t1} and Proposition \ref{t2}, for
   any $y,y'\in Y$, we have
   \begin{align*}
    \dist_N(g_f(y),g_f(y')) \leq \ &\dist_1^W
   ((f^y)_{\ast}(\mu_X),(f^{y'})_{\ast}(\mu_X))\\ 
  = \ & \sup\Big\{\int_{N}\psi(z)  d (f^y)_{\ast}(\mu_X)(z)-\int_N \psi(z) 
    d (f^{y'})_{\ast}(\mu_X)(z) \Big\}\\
    =\ & \sup \Big\{\int_{X} \psi ( f(x,y))  d\mu_X (x) -  \int_{X} \psi (
    f(x,y'))  d\mu_X (x) \Big\}\tag*{} \\
    \leq \ & \int_{X} \dist_N(f(x,y),f(x,y'))  d\mu_X(x), \tag*{}
    \end{align*}where each supremum is taken over all $1$-Lipschitz
   function $\psi:N \to \mathbb{R}$. Observe that the right-hand side of the above
   inequality converges to zero as $\dist_Y(y,y')\to 0$. This completes the proof.
  \end{proof}
  \end{lem}

  \begin{lem}\label{l4}We have $(g_f)_{\ast}(\mu_Y) \in \mathcal{P}_1(N)$.
   \begin{proof}Taking any point $p_0\in N$, from Corollary \ref{c2}, we obtain
    \begin{align*}
     \int_Y \dist_N(\mathbb{E}(f^y), p_0)d\mu_Y(y) \leq \int_{X\times Y}
     \dist_N(f(x,y),p_0)d(\mu_X  \times \mu_Y) (x,y)<+\infty.
     \end{align*}This completes the proof.
    \end{proof}
   \end{lem}
The following example asserts that the equality $\mathbb{E}(f)=
   \mathbb{E}_{y}(\mathbb{E}(f^y))$ does not hold for non-linear CAT(0)-spaces
   in general:

  \begin{ex}\label{e1}\upshape 
For $i=1,2,3$, let $T_i:=\{ (i,r) \mid r\in [0,+\infty)\}$ be a copy of
   $[0,+\infty)$ equipped with the usual Euclidean distance
   function. The \emph{tripod} $T$ is the metric space obtained by gluing together all these
   spaces $T_i$, $i=1,2,3$, at their origins with the intrinsic distance
   function. Let $\{ a,b\}$ be an arbitrary two-point mm-space equipped with the
   uniform probability measure. Let us consider the map $f:\{ a,b\}^2 \to T$ defined by
   $f(a,a):=(1,1)\in T_1$, $f(b,a):=(2,1)\in T_2$, and
   $f(a,b)=f(b,b):=(3,1)\in T_3$. In this case, we easily see that
   $\mathbb{E}(f)=(0,0)$, $\mathbb{E}(f^a)=(0,0)$,
   $\mathbb{E}(f^b)=(3,1)$, and therefore $\mathbb{E}_y(\mathbb{E}(f^y))=(3,1/2)$.
  \end{ex}

\begin{proof}[Proof of Theorem \ref{th1}]
Iterating Corollary \ref{c2}, we have
   \begin{align*}
    \dist_N(\mathbb{E}(f),\mathbb{E}_y(\mathbb{E}(f^y)))\leq \ &
    \int_{X\times Y} \dist_N(f(x,y'),
    \mathbb{E}_y(\mathbb{E}(f^y)))d(\mu_X \times \mu_Y)(x,y')\\
    \leq \ & \int_{X \times Y \times Y}
    \dist_N(f(x,y'),\mathbb{E}(f^{y''}))d(\mu_X \times \mu_Y \times
    \mu_Y)(x,y',y'')\\
    \leq \ & V_1(f).
    \end{align*}Thereby, we obtain the inequality (\ref{s1}).

 To prove the inequality (\ref{s2}), we are
 going to iterate Proposition \ref{p3}. Since $f_{\ast}(\mu_X \times
 \mu_Y)\notin \mathcal{P}_2(N)$ implies $V_2(f)=+\infty$, we assume that
 $f_{\ast}(\mu_X \times \mu_Y)\in \mathcal{P}_2(N)$. From Proposition
 \ref{p3}, we have
 \begin{align}\label{s6}
   &\int_{X\times Y} \dist_N(f(x,y'), \mathbb{E}(f))^2 d(\mu_X \times
   \mu_Y)(x,y')\\
   = \ &\int_Y d\mu_Y(y')\int_X \dist_N(f^{y'}(x), \mathbb{E}(f))^2
   d\mu_X(x) \tag*{}\\
   \geq \ &  \int_Y d\mu_Y(y')\Big\{
   \int_X\dist_N(f^{y'}(x),\mathbb{E}(f^{y'}))^2d\mu_X(x)+
   \dist_N(\mathbb{E}(f^{y'}), \mathbb{E}(f))^2\Big\}\tag*{}\\
   = \ &\int_{X\times Y} \dist_N(f^{y'}(x), \mathbb{E}(f^{y'}))^2d
   (\mu_X \times \mu_Y)(x,y')+\int_Y\dist_N
   (\mathbb{E}(f^{y'}),\mathbb{E}(f))^2 d\mu_Y(y')\tag*{}\\
   \geq \ & \int_{X\times Y} \dist_N(f^{y'}(x), \mathbb{E}(f^{y'}))^2d
   (\mu_X \times \mu_Y)(x,y') +\int_Y\dist_N
   (\mathbb{E}(f^{y'}),\mathbb{E}_y(\mathbb{E}(f^y)))^2
   d\mu_Y(y')\tag*{}\\
   \ & \hspace{10cm}
   +\dist_N(\mathbb{E}(f),\mathbb{E}_y(\mathbb{E}(f^y)))^2\tag*{}.
   \end{align}
 Since
 \begin{align*}
  \dist_N(f^{y'}(x),\mathbb{E}(f^{y'}))^2 +
  \dist_N(\mathbb{E}(f^{y'}),\mathbb{E}_y(\mathbb{E}(f^y)))^2 \geq
  \frac{1}{2}\dist_N(f^{y'}(x), \mathbb{E}_y(\mathbb{E}(f^y)))^2,
  \end{align*}substituting this into
 the inequality (\ref{s6}), we get
 \begin{align*}
  &\int_{X\times Y}\dist_N(f(x,y'),\mathbb{E}(f))^2d(\mu_X \times
  \mu_Y)(x,y')\\ \geq \ &\frac{1}{2}\int_{X\times
  Y}\dist_N(f(x,y'),\mathbb{E}_y(\mathbb{E}(f^y)))^2d(\mu_X \times
  \mu_Y)(x,y') + \dist_N(\mathbb{E}(f), \mathbb{E}_y(\mathbb{E}(f^y)))^2\tag*{}
  \end{align*}
 Since
  \begin{align*}
   &\int_{X\times Y} \dist_N(f(x,y'), \mathbb{E}(f))^2 d(\mu_X \times
   \mu_Y)(x,y')\\ \leq \ &\int_{X\times Y} \dist_N(f(x,y'), \mathbb{E}_y(\mathbb{E}(f^y)))^2 d(\mu_X \times
   \mu_Y)(x,y'),
   \end{align*}we therefore obtain
 \begin{align}\label{s1'}
  \dist_N(\mathbb{E}(f), \mathbb{E}_y(\mathbb{E}(f^y)) )^2 \leq
  \frac{1}{2}\int_{X\times Y} \dist_N(f(x,y'),
  \mathbb{E}_y(\mathbb{E}(f^y)))^2 d(\mu_X \times \mu_Y)(x,y').
  \end{align}

 By virtue of Proposition \ref{p3}, we also get
  \begin{align}\label{s4}
    &\int_{X\times Y} \dist_N (f(x,y'),
   \mathbb{E}_{y}(\mathbb{E}(f^y)))^2  d(\mu_X \times \mu_Y)(x,y')\\
    \leq \ & \int_{X\times Y} d(\mu_X \times \mu_Y)(x,y')\Big\{ \int_Y
   \big\{ \dist_N (f^{y'}(x), \mathbb{E}(f^{y''}))^2-\dist_N
   (\mathbb{E}(f^{y'}), \mathbb{E}_y(\mathbb{E}(f^y)))^2\big\}
   d\mu_Y(y'')\Big\}  \tag*{} \\
   = \ & \int_{X\times Y \times Y} \dist_N(f^{y'}(x),
   \mathbb{E}(f^{y''}))^2 d(\mu_X \times \mu_Y \times \mu_Y) (x,y',y'')
   \tag*{} \\
   & \hspace{8cm} -\int_Y \dist_N(\mathbb{E}(f^{y'}),\mathbb{E}_y(\mathbb{E}(f^y)))^2 d\mu_Y(y').\tag*{} 
   \end{align}Since
   \begin{align*}
     \dist_N(f^{y'}(x), \mathbb{E}(f^{y''}))^2
    \leq \int_X \big\{\dist_N (f^{y''}(x'),f^{y'}(x))^2 -
    \dist_N(f^{y''}(x'), \mathbb{E}(f^{y''}))^2 \big\}d\mu_X(x')
    \end{align*}from Proposition \ref{p3}, substituting this into (\ref{s4}), we have
\begin{align*}
 &\int_{X\times Y} \dist_N (f(x,y'),
   \mathbb{E}_{y}(\mathbb{E}(f^y)))^2  d(\mu_X \times \mu_Y)(x,y')\\ \leq \ &
 V_2(f)^2
 - \int_{X \times Y}\dist_N(f^{y'}(x),
 \mathbb{E}(f^{y'}))^2d(\mu_X \times \mu_Y)(x,y') - \int_Y
 \dist_N(\mathbb{E}(f^{y'}),\mathbb{E}_y(\mathbb{E}(f^y)))^2
 d\mu_Y(y')\tag*{} \\
 \leq \ &  V_2(f)^2 - \frac{1}{2}\int_{X\times Y} \dist_N(f^{y'}(x),
 \mathbb{E}_y(\mathbb{E}(f^y)))^2d(\mu_X \times \mu_Y)(x,y')   \tag*{}.
 \end{align*}We therefore obtain
 \begin{align*}
  \int_{X\times Y}\dist_N(f(x,y'), \mathbb{E}_y
  (\mathbb{E}(f^{y'})))^2 d(\mu_X \times \mu_Y)(x,y') \leq
  \frac{2}{3}V_2(f)^2.
  \end{align*}Combining this with the inequality (\ref{s1'}), we finally obtain the inequality (\ref{s2}). This completes the proof.
 \end{proof}

\section{Applications}
\subsection{Product inequalities}

 \begin{prop}[{Y. G. Reshetnyak, cf.~\cite[Proposition 2.4]{sturm}}]\label{p5}For
  any four points $x_1,x_2, x_3,x_4 $ in a CAT(0)-space $N$, we have
  \begin{align*}\dist_N(x_1,x_3)^2 + \dist_N (x_2,x_4)^2 \leq \dist_N(x_1,x_2)^2
   +\dist_N (x_2,x_3)^2 + \dist_N(x_3,x_4)^2 + \dist_N(x_4,x_1)^2 .
   \end{align*}
  \end{prop}

 Given an mm-space $X$ and a metric space $Y$ we define
\begin{align*}
\obin_Y(X):= \sup \{   V_p(f) \mid f:X\to Y \text{ is a }1
 \text{-Lipschitz map}                  \},
\end{align*}and call it the \emph{observable} $L^p$\emph{-variation} of
 $X$. The idea of the observable $L^p$-variation comes from the quantum and statistical
	mechanics, that is, we think of $\mu_X$ as a state on a configuration
	space $X$ and $f$ is interpreted as an observable.

  \begin{cor}\label{c3}Let $X$ and $Y$ be two mm-spaces and $N$ a
   CAT(0)-space. Then, we have
   \begin{align}\label{s7}
    \obinin_N(X\times Y)^2 \leq \obinin_N(X) + \obinin_N(Y).
    \end{align}
   \begin{proof}Let $f:X \times Y \to N$ be an arbitrary $1$-Lipschitz
    map. Then, putting $Z:=X\times Y$, from Proposition \ref{p5}, we obtain
    \begin{align*}
     V_2(f)^2
     =\ & \frac{1}{2}\int_{Z \times Z} \{
     \dist_{N}(f(x,y),f(x',y'))^2 + \dist_N (f(x,y'), f(x',
     y))^2\}d(\mu_Z \times \mu_Z)(x,y,x',y')\\
     \leq \ & \frac{1}{2}\int_{Z\times Z} \{ \dist_{N}(f(x,y),f(x',y))^2 + \dist_N (f(x',y), f(x',
     y'))^2\\
    & \hspace{1cm}+ \dist_{N}(f(x',y'),f(x,y'))^2 + \dist_N (f(x,y'), f(x,
     y))^2            \}d(\mu_Z \times \mu_Z)(x,y,x',y')\\
     = \ & \int_{X} V_2(f^x)^2d\mu_X(x) +\int_{Y} V_2(f^y)^2d\mu_Y(y)\\
     \leq \ & \obinin_N(X)^2 + \obinin_N(Y)^2.
     \end{align*}This completes the proof.
    \end{proof}
   \end{cor}

  \begin{lem}\label{l5}Let $X$ and $Y$ be two mm-spaces and $Z$ a metric
   space. Then, for any $p\geq 1$, we have
    \begin{align}\label{s8}
     \obin_Z(X\times Y)^p\leq 2^{p-1} \obin_Z(X)^p +2^{p-1} \obin_Z(Y)^p.
     \end{align}
   \begin{proof}Given any $1$-Lipschitz map $f:X\times Y \to Z$,
    putting $W:= X\times Y$, we have
    \begin{align*}
     V_p(f)^p \leq \ &\int_{W\times W}2^{p-1} \{\dist_Z(f(x,y),f(x,y'))^p +
     \dist_Z(f(x,y'), f(x',y'))^p\}d(\mu_W \times \mu_W)(x,y,x',y')\\
     = \ & 2^{p-1}\int_{X}V_p(f^x)^pd\mu_X(x) +
     2^{p-1}\int_{Y}V_p(f^y)^pd\mu_Y(y)\\
     \leq \ & 2^{p-1} \obin_Z(X)^p + 2^{p-1}\obin_Z(Y)^p.
      \end{align*}This completes the proof.
    \end{proof}
    \end{lem}
Note that the inequality (\ref{s7}) is sharper than the
    inequality (\ref{s8}) in the case where $p=2$ and $Z$ is a CAT(0)-space.

    Combining Theorem \ref{th1} and Lemma \ref{l5} we obtain the following corollary:
     \begin{cor}Let $\{ X_n \}_{n=1}^{\infty}$ and $\{
    Y_n\}_{n=1}^{\infty}$ be a sequences of mm-spaces and $\{
    N_n\}_{n=1}^{\infty}$ be a sequences of CAT(0)-spaces. Then, assuming that
    \begin{align*}
     \obvar_{N_n}(X_n) \to 0 \text{ as }n\to \infty \text{ and }
     \obvar_{N_n}(Y_n) \to 0 \text{ as }n\to \infty,
     \end{align*}we have
    \begin{align*}
     \dist_{N_n}(\mathbb{E}(f) , \mathbb{E}_{y_n}(\mathbb{E}(f^{y_n})))
     \to 0 \text{ as }n\to \infty
     \end{align*}for any sequence $\{ f_n :X_n \times
    Y_n \to N_n\}_{n=1}^{\infty}$ of $1$-Lipschitz maps.
    \end{cor}

    \subsection{The non-zero first eigenvalue of Laplacian and the observable $L^2$-variation}
    Although the same method in
    \cite{gromovcat} and \cite{gromov} implies the following
    proposition, we prove it for the completeness. 
  \begin{prop}[{cf.~\cite[Section 13]{gromovcat}, \cite[Section
   $3\frac{1}{2}.41$]{gromov}}]\label{p6}Let $M$ be a compact connected Riemannian manifold
   and $N'$ an $n$-dimensional Hadamard manifold. Then, we have
   \begin{align*}
    \obinin_{N'}(M)\leq 2\sqrt{\frac{n}{\lambda_1(M)}}.
    \end{align*}
   \begin{proof}Let $f:M \to N'$ be an arbitrary $1$-Lipschitz map. We
    shall prove that
    \begin{align}\label{s9}
     \int_M \dist_{N'}(f(x),\mathbb{E}(f))^2d\mu_M (x)\leq \frac{n}{\lambda_1(M)}.
     \end{align}If the inequality (\ref{s9}) holds, then we finish the proof since
    \begin{align*}
     V_2(f)\leq 2\Big( \int_M   \dist_{N'}(f(x),\mathbb{E}(f))^2
     d\mu_M(x)              \Big)^{1/2}\leq 2\sqrt{\frac{n}{\lambda_1(M)}}.
     \end{align*}
    Suppose that
    \begin{align}\label{s10}
     \int_M \dist_{N'}(f(x),\mathbb{E}(f))^2d\mu_M (x)> \frac{n}{\lambda_1(M)}.
     \end{align}We identify the tangent space of $N'$ at
    the point $\mathbb{E}(f)$ with the Euclidean space $\mathbb{R}^n$ and
    consider the map $f_0:= \exp^{-1}_{\mathbb{E}(f)} \circ f:M \to
    \mathbb{R}^n$. According to the hinge theorem (see \cite[Chapter \Roman{yon},
    Remark 2.6]{sakai}), the map $f_0$ is a
    $1$-Lipschitz map. Note that $|f_0 (x)| =
    \dist_{N'}(f(x),\mathbb{E}(f))$ for any $x\in X$ because the map $\exp_{\mathbb{E}(f)}^{-1}$ is
    isometric on rays issuing from the point $\mathbb{E}(f)$. Hence,
    from the inequality (\ref{s10}), we have
    \begin{align*}
     \int_M|f_0(x)|^2d\mu_M (x)> \frac{n}{\lambda_1(M)}.
     \end{align*}Denoting by $(f_0(x))_i$ the $i$-th component
    of $f_0(x)$, we therefore see that there exists $i_0$ such that
    \begin{align}\label{s11}
     \int_M |(f_0(x))_{i_0}|^2d\mu_M(x) > \frac{1}{\lambda_1(M)}.
     \end{align}Note the function $(f_0)_{i_0}$ has the mean zero from
    Lemmas \ref{l1} and \ref{l2}. Combining this with the inequality (\ref{s11}), we therefore obtain
    \begin{align*}
     \lambda_1(M)= \inf \frac{\int_M |\grad_x g|^2d\mu_M(x)}{ \int_M
     |g(x)|^2 d\mu_M(x)}< \lambda_1(M),
     \end{align*}where the infimum is taken over all Lipschitz function
    $g:M \to \mathbb{R}$ with the mean zero. This is a
    contradiction. This completes the proof.
    \end{proof}
   \end{prop}One can obtain a similar result to Proposition \ref{p6} for a finite connected graph.

   \begin{thm}[{cf.~\cite[Proposition 5.7]{funano2}}]\label{t3}Let $X$ be an mm-space and $T$ an
    $\mathbb{R}$-tree. Then, we have
    \begin{align*}
     \obinin_T(X)^2 \leq (38 +16\sqrt{2}) \obinin_{\mathbb{R}}(X)^2. 
     \end{align*}
    \end{thm}
Combining Proposition \ref{p6} with Theorem \ref{t3}, we obtain the following corollary:
   \begin{cor}\label{c4}Let $M$ be a compact connected Riemannian manifold and $T$ an
    $\mathbb{R}$-tree. Then, we have
    \begin{align*}
     \obinin_T(M)^2 \leq \frac{4(38+16\sqrt{2})}{\lambda_1(M)}.
     \end{align*}
    \end{cor}

\begin{proof}[Proof of Corollary \ref{c1}]The corollary follows from
 Theorem \ref{th1}
 together with Corollary \ref{c3}, Proposition \ref{p6}, and Corollary \ref{c4}. This completes the proof.
 \end{proof}
 \begin{ack}\upshape The author would like to express his thanks to
  Professor Takashi Shioya for his valuable comments and assistances during the preparation of
  this paper.
 \end{ack}

	\end{document}